\documentclass[12pt]{amsart}
\usepackage{a4wide,amsmath,amssymb,amsfonts,amscd,verbatim,
psfrag,subfigure,graphicx}

\newcommand{\lebn}

\newtheorem{proposition}[equation]{Proposition}

\newtheorem{theorem}[equation]{Theorem}
\newtheorem{corollary}[equation]{Corollary}
\newtheorem{lemma}[equation]{Lemma}

\theoremstyle{definition}

\newtheorem{definition}[equation]{Definition}

\newtheorem{example}[equation]{Example}

\numberwithin{equation}{section}

\newcommand{\C}{\mathbb{C}}

\newcommand{\R}{\mathbb{R}}
\newcommand{\RE}{\mathbb{R}_{\geq}}
\newcommand{\Z}{\mathbb{Z}}
\newcommand{\SE}{\mathcal{S}}

\newcommand{\PE}{\mathcal{P}}

\newcommand{\B}{\mathcal{B}}
\newcommand{\CE}{\mathcal{C}}

\begin{document}
\bibliographystyle{plain}
\title[Bott towers, crosspolytopes and torus actions \lebn]
{Bott towers, crosspolytopes and torus actions}
\author{Yusuf Civan}
 \address{Department of Mathematics, Suleyman Demirel University,
Isparta, 32260, Turkey}   
\email{ycivan@fef.sdu.edu.tr}


%date{}

\begin{abstract}
We study the geometry and topology of Bott towers in the context of toric
geometry. We show that any $k$th stage of a Bott tower is a smooth 
projective toric variety associated to a fan arising from a crosspolytope; 
conversely, we prove that any toric variety associated to a fan 
obtained from a crosspolytope actually gives rise to a Bott tower. The 
former leads us to a description of the tangent bundle of the $k$th stage 
of the tower, considered as a complex manifold, which splits into a 
sum of complex line bundles. Applying Danilov-Jurkiewicz
theorem, we compute the cohomology ring of any $k$th stage, and
by way of construction, we provide all the monomial identities defining 
the related affine toric varieties. 
\end{abstract}

\maketitle
\section{Introduction}\label{intro}

The theory of toric varieties offers a remarkable area for studying
algebro-geometric and topological problems in the language of combinatorial
objects, called fans, much like simplicial complexes. Many of the algebraic or
topological properties of toric varieties are encoded in the associated fans.
In this direction, our purpose here is to investigate the toric structure of
some complex manifolds, known as Bott towers. These manifolds include such
families of complex manifolds as the Bott-Samelson varieties, and as explained 
by Grossberg and Karshon \cite{GK}, their study combines areas such as
representation theory, combinatorics and algebraic geometry.

Since the theory of toric varieties has been widely studied in currently active 
areas of mathematics, it is impossible to find either a commonly accepted 
definition or a fixed notation. There are two common approaches for constructing 
a toric variety from given combinatorial data. The first  
is due to algebraic geometers in which it is formed by gluing 
affine algebraic varieties; it has the advantage of capturing the 
essence of turning a combinatorial object into an
algebraic one, namely that its local structure resembles the whole.  
The second owes its existence to symplectic geometers, 
and it is more explicit in the sense that the variety can be constructed
as the quotient of a complex space by the action of an algebraic torus. 
Several comprehensive textbooks are now available, by authors such as 
Ewald \cite{GE}, Fulton \cite{Ful}, and Oda \cite{TO}. Throughout, however,
we follow Batyrev  and Ewald, so we refer readers to
\cite{Bat1} and \cite{GE} for background notation.

By a definition, a Bott tower is defined inductively as an iterated
projective bundle so that each stage of the tower is of the form
$\C P(\C \oplus \xi)$ for an arbitrarily chosen line bundle $\xi$ over the
previous stage. It produces a sequence of fibered projective spaces with 
fibers isomorphic to $\C P^1$. A classic example is provided by taking 
$\xi$ to be the trivial line bundle, so each stage of the tower, thus 
obtained, is the product of projective lines. In Section \ref{tower}, 
by incorporating Grossberg and Karshon's construction, we present 
the $k$th stage $N_k$ of any Bott tower of height $n$ as a smooth 
projective toric variety for each $1\leq k\leq n$, describing explicitly 
the associated smooth fan arising from a crossploytope. To achieve this, 
we introduce the notion of \emph{Bott numbers}, which may be of 
particular interest to combinatorialists. Some properties of these 
numbers allow us to reveal the defining Laurent monomials of affine 
toric varieties associated to Bott towers. We also provide an example 
showing that Bott towers are not Fano varieties in general. We note that a
Bott tower of height $2$ is actually a Hirzebruch surface.
Moreover, we show that bounded flag manifolds are also examples of 
Bott towers, and the results of~\cite{BuR1} therefore suggest that 
they might have a role to play in complex bordism and cobordism theory 
that has yet to be revealed.

Toric geometry is something of a two-way study in the sense that we may start with a
normal algebraic variety which contains the algebraic torus as a dense open 
subset, and then recover the associated combinatorial data. Conversely, 
beginning with a fan, we can construct such a variety. Following the latter
pattern, we prove that any toric variety associated to a smooth fan arising 
from a crosspolytope is actually a Bott tower. 

One of the advantage of having the toric structure of Bott towers is 
that we may easily describe their tangent bundles. According to Ehlers \cite{Ehl}, 
the tangent bundle of a toric variety splits as a sum of line bundles 
(known as the \emph{hyperplane bundles}), and each of these bundles is 
obtained as the pull-back of the normal bundle of a codimension-one 
subvariety corresponding to a one-dimensional cone of the associated fan.  
In a parallel work \cite{YC3}, we use such a splitting to assist our
computation of real and complex $K$-groups of Bott towers and determine
various structures (almost and stably complex structure, etc.) on them.

\textit{Acknowledgments}. I would like to thank Professor Nige Ray for 
many useful discussions and for his generous encouragement.
%%%%%%%%%%%%%%%%%%%%%%%%%%%%%%%%%%%%%%%%%%%%%%%%%%%%%%%%%%%%%%%%%%%
%%%%%%%%%%%%%%%%%%%%%%%%%%%%%%%%%%%%%%%%%%%%%%%%%%%%%%%%%%%%%%%%%%%
\section{Bott Towers}\label{tower}

We begin with a precise description of Bott towers. The
way we display these objects fits into the setting of~\cite{Ray}, and we
provide an alternative proof for Grossberg and Karshon's 
construction~\cite{GK}.  

\begin{figure}[ht]
\begin{picture}(100, 150)(-1, 17)
\put(108, 132){\makebox(0,0){$\C P(\C\oplus \xi_{n-1})= N_n$}}
\put(88, 122){\vector(-1, -1){14}}
\put(92, 112){\makebox(0,0){$\pi_n$}}
\put(72, 104){\makebox(0,0){$\cdot$}}
\put(70, 102){\makebox(0,0){$\cdot$}}
\put(68, 100){\makebox(0,0){$\cdot$}}
%%%
\put(65, 97){\vector(-1, -1){14}}
\put(68, 89){\makebox(0,0){$\pi_3$}}
%%%
\put(52, 75){\makebox(0,0){$\C P(\C\oplus \xi_1)= N_2$}}
%%%
\put(38, 66){\vector(-1, -1){14}}
\put(43, 56){\makebox(0,0){$\pi_2$}}
%%%%%%
\put(28, 42){\makebox(0,0){$\C P(\C\oplus \xi_0)= N_1$}}
%%%%
\put(9, 31){\vector(-1, -1){14}}
\put(16, 22){\makebox(0,0){$\pi_1$}}
%%%%%
\put(4, 9){\makebox(0,0){$\ast= N_0$}}
\end{picture}
\newline
\caption{A Bott tower of height $n$}\label{BT}
\end{figure}
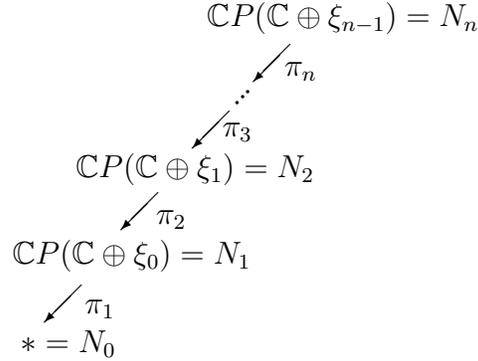

Let $\xi_0$ be a trivial line bundle over a single point $N_0= *$, and let 
$N_1:=\C P(\C\oplus \xi_0)=\C P^1$. Similarly, we may choose
any holomorphic line bundle $\xi_1$ over $N_1$, and 
take its direct sum with the trivial line bundle. By projectifying each
fiber, we obtain a manifold $N_2=\C P(\C\oplus \xi_1)$, 
which is a bundle over
$N_1$ with fiber $\C P^1$. We may repeat this process $n$ times, so
that each $N_k$ is a $\C P^1$-bundle over $N_{k-1}$ (see Figure \ref{BT}).
We can consider each $N_k$ as the space of lines in 
$\C\oplus \xi_{k-1}$.
If we think of $\C P^1$ as a sphere with a south 
pole $[1, 0]$ and a north pole $[0, 1]$, 
then the zero section of $\xi_{k-1}$ gives rise to a 
holomorphic section, the  \emph{south pole section}:
\begin{displaymath} 
i^S_k \colon N_{k-1}\rightarrow N_k;
\end{displaymath}
similarly, we may obtain the \emph{north pole section} 
$i^N_k \colon N_{k-1}\rightarrow N_k$
by allowing the first coordinate in 
$\C P(\C\oplus \xi_{k-1})$ to vanish.
\begin{definition}
A \emph{Bott tower of height $n$} is defined to be  a collection of complex
manifolds $\{N_k\colon k\leq n\}$, constructed by the above process. The
bundles $\xi_0, \ldots, \xi_{n-1}$ are called the \emph{associated line
bundles} of the tower.   
\end{definition}
\begin{example}
The collection 
$\{(\C P^1)^k:=\C P^1\times \ldots \times \C P^1 
\textrm{($k$ times)}\colon k\leq n\}$ 
is a Bott tower. In this case, each associated line bundle
$\xi_k$ is chosen to be trivial for every $0\leq k\leq n-1$, and $\pi_k$ is the
obvious projection.
\end{example}
Our next example has a particular importance in complex cobordism theory and the 
resulting manifolds are known as  bounded flag manifolds (see~\cite{BuR1}). 
They were, in fact, introduced by Bott $\&$ Samelson \cite{BS}.

We write $[n]$ for the set of natural numbers $\{1, 2, \ldots, n\}$ with the 
standard linear ordering and an interval in the poset $[n]$ has the form
$[i, j]$ for some $1\leq i\leq j\leq n$. Throughout, $\omega_1, \ldots, \omega_{n+1}$ 
will denote the standard basis vectors in $\C^{n+1}$, and we write $\C_I$ for the subspace 
spanned by the vectors $\{\omega_i\colon i\in I\}$, where $I\subset [n+1]$.  
\begin{definition}
A flag $U\colon 0< U_1<\ldots < U_n <\C^{n+1}$ is called \emph{bounded} if
$\C_{[i-1]}< U_i$ for each $1\leq i\leq n$. The space of all bounded flags in
$\C^{n+1}$ is called \emph{bounded flag manifold}, which is an 
$n$-dimensional smooth complex manifold and will
be denoted by $B(\C^{n+1})$ (or simply by $B_n$ when there is no confusion).
\end{definition}  
As a consequence of the definition, each factor $U_i$ of any bounded flag 
$U\in B(\C^{n+1})$ is of the form $\C_{[i-1]}\oplus L_i$, where $L_i$ is a
line in $\C_i\oplus L_{i+1}$ for $1\leq i\leq n$, and
$L_{n+1}=\C_{n+1}$. Therefore, we may display $U$ as
\begin{equation}\label{EQ6}
U\colon 0 < L_1 <\C_1\oplus L_2< \ldots <\C_{[n-1]}\oplus L_n < \C^{n+1}.
\end{equation}
We define maps $q_i$ and $r_i\colon B(\C^{n+1})\rightarrow \C P_{[i, n+1]}$ by
letting $q_i(U)=L_i$ and $r_i(U)=L_i^{\bot}$, where $L_i^{\bot}$ is the
orthogonal complement of $L_i$ in   
$\C_i\oplus L_{i+1}$ for each $U\in B(\C^{n+1})$, and  
$1\leq i\leq n$. We consider complex line bundles $\eta_i$ and 
$\eta_i^{\bot}$ over $B_n$,  classified respectively by the maps 
$q_{n-i+1}$ and $r_{n-i+1}$ for every  $1\leq i\leq n$, and we set 
$\eta_0$ to be the trivial line bundle with fiber $\C_{n+1}$. We 
sometimes refer them as the \emph{canonical line bundles} on $B_n$.
\begin{theorem}\label{T1}
The collection $\{B(\C_{[n-k+1, n+1]})\colon k\leq n\}$ of bounded flag
manifolds is a Bott tower of height $n$, 
where $B(\C_{[n-k+1, n+1]})$ denotes the set of bounded flags in 
$\C_{[n-k+1, n+1]}$. 
\end{theorem}
\begin{proof}
As above, we may define maps 
\begin{equation*}
q_i\quad\mathrm{and}\quad r_i\colon B(\C_{[n-k+1, n+1]})\rightarrow 
\C P_{[i, n+1]}\quad \mathrm{for}\; n-k+1\leq i\leq n+1,
\end{equation*} 
and the respective complex line bundles $\eta_0, \eta_1, \ldots, \eta_k$ 
over $B(\C_{[n-k+1, n+1]})$.
We proved in~\cite{YC1} that the complex manifolds 
$\C P(\C_{n-k+1}\oplus \eta_{k-1})$ 
and  $B(\C_{[n-k+1, n+1]})$ are diffeomorphic
for each $1\leq k\leq n$. Therefore, we can take
$B(\C_{[n-k+1, n+1]})$ to be the $k$th stage of the tower; hence, we may  
abbreviate $B(\C_{[n-k+1, n+1]})$ to $B_k$. The projection 
$\pi_k\colon B_k\rightarrow B_{k-1}$ maps each flag  
\begin{align*}
U_k\colon &0< L_{n-k+1}<\C_{n-k+1}\oplus L_{n-k+2}<\ldots
<\C_{[n-k+1, n-1]}\oplus L_n < \C_{[n-k+1, n+1]}\\
\intertext{in $B_k$ to the flag}
U_{k-1}\colon &0< L_{n-k+2}<\C_{n-k+2}\oplus L_{n-k+3}<\ldots
<\C_{[n-k+2, n-1]}\oplus L_n < \C_{[n-k+2, n+1]}
\end{align*} 
in $B_{k-1}$, whose fiber consists of the lines in $\C_{n-k+1}\oplus
L_{n-k+2}$, and is therefore isomorphic to $\C P^1$. 
The south and north pole sections $i^S_k$ and 
$i^N_k\colon B_{k-1}\rightarrow B_k$ are given respectively by
\begin{align*}
i^S_k(U_{k-1}):&=0<\C_{n-k+1}< \C_{n-k+1}\oplus L_{n-k+2}<\ldots
<\C_{[n-k+1, n-1]}\oplus L_n < \C_{[n-k+1, n+1]},
\intertext{and}
i^N_k(U_{k-1}):&=0< L_{n-k+2}<\C_{n-k+1}\oplus L_{n-k+2}<\ldots
<\C_{[n-k+1, n-1]}\oplus L_n < \C_{[n-k+1, n+1]}.
\end{align*} 
\end{proof}
%%%%%%%%%%%%%%%%%%%%%%%%%%%%%%%%%%%%%%%%%%%%%%%%%%%%%%%%%%%%%%%%%%%%%%%%%
%%%%%%%%%%%%%%%%%%%%%%%%%%%%%%%%%%%%%%%%%%%%%%%%%%%%%%%%%%%%%%%%%%%%%%%%%
\section{Bott towers as toric varieties}
As we explained earlier, one way of constructing toric varieties is to display
them as the quotient of a complex space by the action of an algebraic torus
(see \cite{Bat1}). To be more explicit, if $\Sigma$ is a smooth and complete 
fan in $\R^n$ with the generating set $G(\Sigma)=\{x_1, \ldots, x_m\}$, 
then each primitive collection $\PE$ in $G(\Sigma)$ defines an affine 
subspace in $\C^m$. Running over all primitive collections, we let 
$\mathcal{U}(\Sigma)$ be the complement of the union of corresponding 
affine subspaces. On the other hand, from the kernel $\mathcal{R}(\Sigma)$ 
of the map $\pi_{\Sigma}\colon \Z^m\rightarrow \Z^m$ defined by 
$\pi_{\Sigma}(e_i):=x_i$, where $e_1, \ldots, e_m$ are the standard basis
vectors of $\R^m$, we obtain an algebraic subtorus $\mathcal{D}(\Sigma)$
that acts on $\mathcal{U}(\Sigma)$, so that the quotient
$\mathcal{U}(\Sigma)/\mathcal{D}(\Sigma)$ is the associated toric
variety $X_{\Sigma}$. 

In the case of Bott towers, Grossberg and Karshon~\cite{GK} have already been
able to describe explicitly how to construct Bott towers as quotients. 
By incorporating their work, we may readily present
any Bott tower as a sequence of smooth projective toric varieties.
We first introduce some notation.

%%%%%%%%%%%%%%%%%%%%%%%%%%%%%%%%%%%%%%%%%%%%%%%%%%%%%%%%%%%%%%%%%%%%%%%%
\subsection{Bott Numbers}
Let $\{c(i, j) \colon 1\leq i< j\leq n\}$ be a collection of $n(n-1)/2$ 
arbitrary integers. Then, by setting $c(i,i):=1$ for each $1\leq k\leq n$, 
we denote the $k$-tuple $(c(1, k), \ldots , c(k, k))$ in 
$\Z^k$ by $c_k$ for each $1\leq k\leq n$. Moreover, we let 
$\mathbf{c}:=(c_1, \ldots, c_n)$ be the resulting integral sequence.

We consider an upper triangular $n\times n$-matrix $\CE(n)$, whose
$k$th row is given by the vector 
$(0,\ldots, 0, c(k,k), c(k, k+1),\ldots, c(k,n))$ for each $1\leq k\leq n$.
In otherwords, if $\CE(n)=(c_{ij})$, then
\begin{equation*}
c_{ij}:=\begin{cases}
c(i, j),& \textrm{if\;$i\leq j$}\\
0, & \textrm{otherwise}.
\end{cases}
\end{equation*}
\begin{definition}
Let $\mathbf{c}=(c_1, \ldots, c_n)$ be an integral sequence. We call the 
upper triangular matrix $\mathcal{B}(n)$ satisfying 
\begin{equation}
\CE(n)^{-1}=-\B(n)
\end{equation} 
the \emph{Bott matrix} associated to $\mathbf{c}$.
Inparticular, if we write $\B(n)=(b(i, j))$, we then have
\begin{align}\label{P1}
\notag & b(i,i)=-1\quad \textrm{and}\quad  b(i,j)=0\;\textrm{if\;$i>j$} \\
& b(i,j)=- \sum_{i< k\leq j}c(i,k)b(k, j)\quad \textrm{if\;$i<j$}.
\end{align}
The numbers $b(i, j)$ for all $1\leq i< j\leq n$ are said to be
the \emph{Bott numbers} associated to $\mathbf{c}$.
\end{definition}

We may readily extend the definition of a Bott matrix to any subset of $[n]$ 
as follows. Let $I=\{i_1, \ldots, i_m\}\subset [n]$ be given
such that $1\leq i_1< \ldots < i_m\leq n$. We then form an upper 
triangular $m\times m$-matrix $\CE(I):=(c_{rs})$ by
\begin{equation*}
c_{rs}:=\begin{cases}
c(i_r, i_s),& \textrm{if\;$r\leq s$}\\
0, & \textrm{otherwise},
\end{cases}
\end{equation*}
and obtain the associated \emph{Bott matrix} $\B(I)$ of $I$ satisfying  
$\CE(I)^{-1}=-\B(I)$. In this case, if we denote the top-right-corner entry 
of $\B(I)$ by $b(I)$, we then call $b(I)$ the \emph{Bott number} of $I$
associated to $\mathbf{c}$. 

There is also an alternative way to describe the Bott number $b(I)$ as follows. 
Let $J\subset [n]$ be given such that $\textrm{min}(J)=i$ and 
$\textrm{max}(J)=j$. Then, for any subset $L=\{l_1,\ldots, l_k\}\subseteq J$
with $l_1<\ldots < l_k$, we associate an integer with it by
\begin{equation}\label{CH2}
p(L,J):=c(i, l_1)c(l_1, l_2)\ldots c(l_k, j),
\end{equation}  
and if $L=\emptyset$, we set $p(\emptyset,J):=c(i, j)$. Moreover, 
if we define $c(J):=\sum_{L\subseteq J}(-1)^{|L|}p(L,J)$, we may easily 
deduce that
\begin{equation}\label{M1}
b(I)=\sum_{J\subseteq I}c(J).
\end{equation} 
for any $I\subset [n]$ such that $|I|> 1$.

We now discuss some properties of Bott numbers associated to some
specific subsets of $[n]$, which we need in Section \ref{LS}. 
Let $\SE^k$ denote the set of binary codes of length $k$, where
$\SE=\{0, 1\}$. For a given $w=w_1\ldots w_k\in \SE^k$,
we define
\begin{equation}
\mathbb{I}^j_i:=\{l\in (i, j)\colon w_l=1\}\quad\mathrm{and}
\quad \mathbb{O}^j_i:=\{l\in (i, j)\colon w_l=0\}
\end{equation}   
for any $1\leq i< j\leq k$. We note that since $\mathbb{O}^j_i\subset [n]$, 
the integer $b(\mathbb{O}^j_i)$ is defined as above.
\begin{lemma}\label{L2}
Let $w\in \SE^k$ be given. Then, for any $1\leq i<j\leq k$, we have 
\begin{align}
\mathbf{(i)}\quad &\sum_{l\in \mathbb{I}_i^j\cup \{i\}} b(\mathbb{O}^j_l) b(i, l) 
+ b(i, j)=0,\label{B2}\\ 
\mathbf{(ii)}\quad &\sum_{l\in \mathbb{I}_i^j\cup \{j\}} b(\mathbb{O}^l_i) b(l, j) 
+ b(i, j)=0.\label{B3}
\end{align}
\end{lemma}
\begin{proof}
It is obvious that proving \eqref{B2} is equivalent to show that
\begin{equation}\label{BE1}
b(i, j)= b(\mathbb{O}^j_i)-\sum_{l\in \mathbb{I}_i^j}
b(\mathbb{O}^j_l) b(i, l)
\end{equation}
which follows from \eqref{M1}.
\end{proof}
%%%%%%%%%%%%%%%%%%%%%%%%%%%%%%%%%%%%%%%%%%%%%%%%%%%%%%%%%%%%%%%%%%%%%%%%%
\subsection{Toric Structures}\label{TS}
Once we have described Bott numbers, we may easily provide the toric structure of
Bott towers. We begin with recalling the definition of our central 
combinatorial objects, namely crosspolytopes (see \cite{Bro}).
\begin{definition}
Let $P^1$ be a line segment in $\R^n$. We proceed by induction, and
assume that $P^k$ is defined for some $k>1$. Let $I_{k+1}$ be a line
segment such that the intersection
\begin{displaymath}
(\textrm{relint}\;P^k )\cap (\textrm{relint}\;I_{k+1})=\{q\},
\end{displaymath}
is a single point. Then we define $P^{k+1}:=\textrm{conv}(P^k\cup I_{k+1})$.
The polytope thus obtained is called a $(k+1)$-\emph{crosspolytope}. It easily
follows that the vertex set of an $n$-crosspolytope $P^n$ can be given as
\begin{displaymath}
V(P^n)=\{x_1, \ldots, x_n, \widehat{x}_1, \ldots, \widehat{x}_n\}
\end{displaymath}
such that the line segment $[x_i, \widehat{x}_i]$ is not an edge of $P^n$ for
all $i=1, \ldots, n$.
\end{definition}

\begin{figure}[ht]
\begin{center}
\psfrag{a}{$(1, 0)$}
\psfrag{b}{$(-1, c(1, 2))$}
\psfrag{c}{$(0, 1)$}
\psfrag{d}{$(0, -1)$}
\includegraphics[width=4in,height=3.6in]{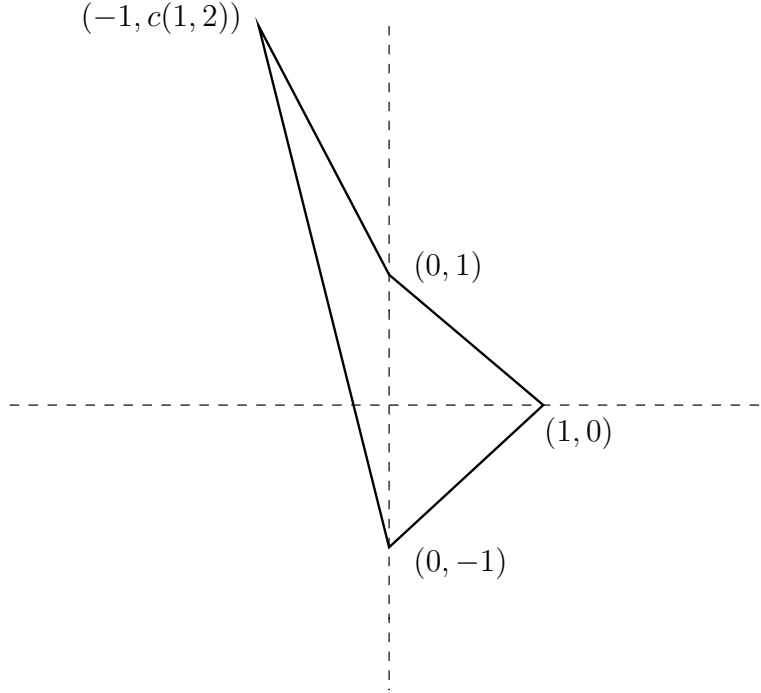}
\end{center}
\caption{The fan of a Bott tower of height $2$, which is not a Fano variety.}\label{FB}
\end{figure}

For a given integral sequence $\mathbf{c}=(c_1, \ldots, c_n)$,
we consider a $k$-dimensional fan $\Sigma(k)$ in $\R^k$ for each 
$k\leq n$ with the generating set 
$G(k)=\{a_{1, 0}, a_{1, 1}, \ldots, a_{k, 0}, a_{k, 1}\}$ defined by
\begin{equation}\label{A2}
a_{i, \gamma}:=\begin{cases}
e_i, & \text{if $\gamma=0$},\\
(0, \ldots, 0, b(i, i), b(i, i+1), \ldots, b(i, k)), & \text{if $\gamma=1$}
\end{cases}
\end{equation}
for any $1\leq i\leq k$. The $k$-dimensional cones of $\Sigma(k)$ are given by
\begin{equation}\label{E2}
\sigma(w):= \RE a_{1, w_1} +\ldots + \RE a_{k, w_k}
\end{equation}
for any $w:=w_1\ldots w_k\in \SE^k$. 
\begin{proposition}
The fan $\Sigma(k)$ is smooth and projective for the integral sequence 
$(c_1, \ldots, c_k)$, which is the first $k$-tuple of $\mathbf{c}$.
\end{proposition}
\begin{proof}
We consider the intersection points of a unit sphere centered 
at the origin with the one-dimensional cones of $\Sigma(k)$. The convex hull of 
those points is clearly a crosspolytope which spans $\Sigma(k)$. Smoothness 
follows from the construction. 
\end{proof}
We note that the smooth toric variety $X_{\Sigma(k)}$ associated to the fan 
$\Sigma(k)$ is not a Fano variety in general (see Figure \ref{FB}).
We may now apply Batyrev construction's (see [\cite{Bat1}, p.4] to the fan $\Sigma(k)$.
Since $\Sigma(k)$ is spanned by a crosspolytope, the space 
$\mathcal{U}(\Sigma(k))$ is $(\C^2\backslash 0)^k$. Moreover, the 
group $\mathcal{R}(\Sigma(k))$ being the kernel of the map 
\begin{align*}
\pi_k \colon &\Z^{2k}\rightarrow \Z^k,\\
&e_{2r-1}\mapsto e_r,\\
&e_{2r} \mapsto (0,\ldots, 0, b(r, r), b(r, r+1), \ldots, b(r, k))
\end{align*}
for $1\leq r\leq k$, is generated by the following $k$ vectors
\begin{align*}
&(1, c(1, 1), 0, c(1, 2), 0, c(1, 3), 0, \ldots, 0, c(1, k)),\\
&(0, 0, 1, c(2,2), 0, c(2, 3), 0, \ldots, 0, c(2, k)),\\
&\vdots\\
&(0, 0, \ldots, 0, 1, c(k-1, k-1), 0, c(k-1, k)),\\
&(0, 0, \ldots, 0, 1, c(k, k)).
\end{align*}
Therefore, the corresponding $k$-dimensional  algebraic torus 
$\mathcal{D}(\Sigma(k))$ is given by
\begin{equation*}
\{(t_1, t_1; t_2, t^{c(1, 2)}_1 t_2; \ldots ; t_k, t^{c(1, k)}_1 t^{c(2, k)}_2
\ldots t^{c(k-1, k)}_{k-1} t_k)\colon t_i\in \C_{\times}\;\textrm{for}\;
1\leq i\leq k\},
\end{equation*}
where $\C_{\times}$ is the algebraic torus. Furthermore the group 
$\mathcal{D}(\Sigma(k))$ acts on $(\C^2\backslash 0)^k$ diagonally.
We denote the resulting smooth and projective toric variety associated
to $\Sigma(k)$ by $X_{\Sigma(k)}$. 

\begin{example}
We consider an integral sequence $\mathbf{c}=(c_1, \ldots, c_n)$ given 
for each $1\leq k\leq n$ by $c_k=(0, 0, \ldots, 0, 1)\in \Z^k$, and let 
$\{M_k\colon k\leq n\}$ be the collection of associated complex manifolds. The
fan $\Sigma(k)$ constructed from $(c_1, \ldots, c_k)$ has the generating set
$G(\Sigma(k))=\{e_1, \ldots, e_k, -e_1, \ldots, -e_k\}$; hence, the associated
toric variety is $(\C P^1)^k$. Therefore, we have $X_{\Sigma(k)}= (\C P^1)^k$ for
any $1\leq k\leq n$.
\end{example}
We next provide an alternative proof of Grossberg and Karshon's construction
in a toric setting. Our method basically applies to the relation between 
linear support functions of fans and holomorphic line bundles over toric
varieties (see \cite{TO}).
\begin{theorem}\label{T2}
For a given integral sequence $\mathbf{c}=(c_1, \ldots, c_n)$, the family
of manifolds $\{X_{\Sigma(k)} \colon k\leq n\}$ is a Bott tower, and any Bott tower arises
in this way.
\end{theorem} 
\begin{proof}
Let $\{X_{\Sigma(k)} \colon k\leq n\}$ be given. We
proceed by induction on $k$. When $k=1$, there is nothing to prove, since 
$X_{\Sigma(1)}= \C P^1$. Therefore, we may assume that there exist line bundles 
$\xi_0, \ldots, \xi_{k-1}$ for some $k< n$ such that 
$X_{\Sigma(l)}=\C P(\C\oplus \xi_{l-1})$ for any $1\leq l\leq k$. 
Then, it can be easily verified that the function 
$h_k\colon |\Sigma(k)|\rightarrow \R$ given by
\begin{equation}\label{A3}
h_k(a_{i, \gamma}):=\begin{cases}
0, & \text{if\;$\gamma=0$},\\
b(i, k+1), & \text{if\;$\gamma=1$}
\end{cases}
\end{equation}
for each $1\leq i\leq k$ is a $\Sigma(k)$-linear support function.
However, each such function $h_k$ defines an equivariant line bundle 
$L(h_k)$ over $X_{\Sigma(k)}$. On the other hand, the space 
$\C P(\C \oplus L(h_k))$ is a smooth projective toric variety whose 
associated fan may be described as the join 
$\Sigma \cdot \widetilde{\Sigma}$, where 
$\Sigma:=\{\RE e_{k+1}, \RE (-e_{k+1}), \{0\}\}$ and  
$\widetilde{\Sigma}:=\Psi(\Sigma)$, and the linear map 
$\Psi \colon \R^k\rightarrow \R^{k+1}$ is given by $y\mapsto (y, h_k(y))$.
Since we have $\Sigma(k+1)=\Sigma \cdot \widetilde{\Sigma}$, we deduce 
that $X_{\Sigma(k+1)}\cong \C P(\C \oplus L(h_k))$, which completes the first
part of the assertion.

Let a Bott tower $\{N_k \colon k\leq n\}$ of height $n$ be given with the
associated line bundles $\xi_0, \ldots, \xi_{n-1}$.
Once again, we proceed by induction on $k$. Since the base case is obvious, we
may assume that there exists an integral sequence $(c_1, \ldots, c_k)$ for
which $N_k=\C P(\C \oplus \xi_{k-1})$ is of the form $X_{\Sigma(k)}$ for some 
$1< k< n$. Let $\xi_k$ be any holomorphic line bundle over $N_k=X_{\Sigma(k)}$, 
and without loss of generality, we assume that it is nontrivial. 
Then, there exists a $\Sigma(k)$-linear support function 
$h\colon |\Sigma(k)|\rightarrow \R$ such that the bundles $L(h)$ and $\xi_k$
are isomorphic. Since $h$ is  $\Sigma(k)$-linear, there exists
$r(w)=(r_1, \ldots, r_k)\in \Z^k$ for any $w=w_1\ldots w_k\in \SE^k$ such that
\begin{equation}\label{E3}
h(x)=\langle r(w), x\rangle
\end{equation} 
for all $x\in \sigma(w)$. Specifically, we consider the binary code 
$w(1):=11\ldots 1$ and its associated integral vector 
$r(1):=(r^1_1, \ldots, r^1_k)\in \Z^k$ satisfying \eqref{E3}.
On the other hand, if $w_i=0$ for any given $w\in \SE^k$ and for some
$1\leq i\leq k$, then $h(a_{i, 0})=\langle r(w), e_i\rangle = r_i=:X_i\in \Z$, 
where $r(w)=(r_1, \ldots, r_k)$ is the associated integral vector of $w$.   
If we define 
\begin{equation*}
c(1, k+1):=X_1-r^1_1, \ldots, c(k, k+1):=X_k-r^1_k, 
\end{equation*}
then we necessarily have 
\begin{equation*}
h(a_{i, \gamma})=\begin{cases}
X_i,& \text{if\;$\gamma=0$},\\
b(i, k+1)-X_i+X_{i+1} b(i, i+1)+\ldots+ X_k b(i, k),& \text{if\;$\gamma=1$}
\end{cases}
\end{equation*}
for each $1\leq i\leq k$ by \eqref{P1}.
Let $\Psi_k\colon \R^k\rightarrow \R^{k+1}$ be a linear map given by
$y\mapsto (y, h(y))$. We define 
\begin{equation*}
\widetilde{\Sigma}:=\{\Psi_k(\sigma)\colon \sigma\in \Sigma(k)\}\quad 
\mathrm{and}\quad \Sigma':=\widetilde{\Sigma}\cdot \Sigma_1,
\end{equation*}
where $\Sigma_1:=\{\RE e_{k+1}, \RE (-e_{k+1}), \{0\}\}$. Then, any
$(k+1)$-dimensional cone of $\Sigma'$ is of the form
\begin{equation*}
\sigma(w)=\RE A_{1, w_1}+\ldots+\RE A_{k+1, w_{k+1}},
\end{equation*}
where 
\begin{equation*}
A_{i, \gamma}:=\begin{cases}
(0, \ldots, 0, 1, 0, \ldots, X_i), & \text{if\;$\gamma=0$}\\
(0, \ldots, 0, -1, f(i, i+1), \ldots, f(i, k), h(a_{i, 1})),& \text{if\;
$\gamma=1$}.
\end{cases}
\end{equation*}
Let $\Sigma(k+1)$ be the fan associated to the integral sequence
$(c_1, \ldots, c_{k+1})$ constructed as in \eqref{A2}. We consider the
following map
\begin{align*}
L \colon &\R^{k+1}\rightarrow \R^{k+1},\\
&e_i\mapsto (0, \ldots, 0, 1, 0, \ldots, 0, X_i)\quad \mathrm{for}\;1\leq
i\leq k,\\
&e_{k+1}\mapsto e_{k+1},
\end{align*}
which is represented by the matrix
\begin{equation*}
\mathbf{A}=
\left(
\begin{matrix}
1       &   0       &   \dots  & 0       & 0\\
0       &   1       &   \dots  & 0       & 0\\
\vdots  &   \vdots  &  \vdots  & \ddots  & \vdots\\
0       &   0       &   \dots  & 1       & 0\\
X_1     &  X_2      &   \dots  & X_k     & 1\\
\end{matrix}
\right).
\end{equation*}
It easily follows that $L$ is unimodular, and it maps the cones of 
$\Sigma(k+1)$ bijectively onto the cones of $\Sigma'$; hence, 
$\Sigma(k+1)\cong \Sigma'$. However, this guarantees that the space
$\C P(\C \oplus L(h))\cong \C P(\C \oplus \xi_k)$ is a toric variety
associated to the fan $\Sigma(k+1)$ with a generating set
$\{a_{1, 0}, a_{1, 1}, \ldots, a_{k+1, 0}, a_{k+1, 1}\}$ defined by
\eqref{A2} for the integral sequence $(c_1, \ldots, c_{k+1})$. 
\end{proof}
We note that every finite poset defines a Bott tower in the following way.
Let $P$ be a finite poset with $n$ elements and let $Z(P)$ be its Zeta matrix
(see \cite{RS}). By taking $\CE(n):=Z(P)$, then the corresponding Bott matrix 
is given by $\B(n)=-M(P)$, where $M(P)$ is the Moebius matrix of $P$, 
from which we can construct a Bott tower of height $n$.
\begin{corollary}\label{C1}
Let $\beta$ be a holomorphic line bundle over the $k$th stage of a Bott tower
$\{N_k\colon k\leq n\}$. Then there exist integers $b_1$, $\ldots$, $b_k$
such that
\begin{equation}
\beta\cong (\C^2\backslash 0)^k\times_{\C^k_{\times}} \C,
\end{equation}
where the action of $\C^k_{\times}$ is given by
\begin{equation}
((x, y), v)\cdot t:=((x, y)\cdot t, t^{b_1}_1\ldots t^{b_k}_k v)
\end{equation}
for any $t=(t_1, \ldots, t_k)\in \C^k_{\times}$.
\end{corollary}
One of the important consequence of Corollary \ref{C1} is that each bundle
$\xi_k$ associated to the $(k+1)$th stage $N_{k+1}$ of the tower 
for $1\leq k\leq n-1$ is given over $N_k$ by the action 
\begin{equation}\label{E4}
((x, y), v)\cdot t:=((x, y)\cdot t, t^{c(1, k+1)}_1\ldots 
t^{c(k, k+1)}_k v)
\end{equation} 
for all $t=(t_1, \ldots, t_k)\in \C^k_{\times}$.
\begin{corollary}
The integral sequence $\mathbf{c}=(c_1, \ldots, c_n)$ associated to Bott tower of bounded
flag manifolds $\{B_k\colon k\leq n\}$ is given by 
$c_k:=(0, \ldots, 0, -1, 1)\in \Z^k$ for each $1\leq k\leq n$.
\end{corollary}
\begin{proof}
In order to retain the notation that was already introduced, 
a slight modification of the quotient description of the varieties 
$X_{\Sigma(k)}$ is required to construct bounded flag manifolds 
as toric varieties. The obvious reason is that the first stage of the tower 
$\{B(\C_{[n-k+1, n+1]})\colon k\leq n\}$ is the projectivization of $\C_n\oplus
\C_{n+1}$. We therefore define $X_{\Sigma(k)}$ for each $k\leq n$ to be the quotient of
$(\C^2\backslash 0)^k$ by the $k$-fold algebraic torus $\C^k_{\times}$, under
the action
\begin{align}\label{A4}
(x_1, y_1;&\ldots;x_k, y_k)\cdot (t_1, \ldots, t_k):=\\\notag
&(x_1 t_1, y_1 t_1 t^{-1}_2; \ldots; x_{k-1} t_{k-1}, y_{k-1} t_{k-1}
t^{-1}_k; x_k t_k, y_k t_k).
\end{align}   
Therefore, the integral sequence $\mathbf{c}=(c_1, \ldots, c_n)$ associated to 
$\{X_{\Sigma(k)}\colon k\leq n\}$ is given by $c_k:=(0, \ldots, 0, -1, 1)\in \Z^k$ for each
$1\leq k\leq n$. 

For a given vector $(x, y)\in (\C^2\backslash 0)^k$, we define
\begin{equation}\label{E5}
\begin{split}
l_{n+1}:&= \omega_{n+1},\\
l_n:&= x_k \omega_n + y_k l_{n+1},\\
&\vdots\\
l_{n-k+1}:&= x_1 \omega_{n-k+1} + y_1 l_{n-k+2}.\\
\end{split}
\end{equation}
If $L_j$ denotes the line in $\C_{[n-k+1, n+1]}$ spanned by 
the vector $l_j$ for $n-k+1\leq j \leq n$, then 
\begin{displaymath}
U(x, y)\colon 0< L_{n-k+1}< \C_{n-k+1}\oplus L_{n-k+2} <\ldots 
<\C_{[n-k+1, n-1]}\oplus L_n < \C_{[n-k+1, n+1]}
\end{displaymath}
is a bounded flag in $B(\C{[n-k+1, n+1]})$ determined by the vector $(x, y)$.
Conversely, for any bounded flag $U\in B(\C{[n-k+1, n+1]})$, we can find
$(x, y)\in (\C^2\backslash 0)^k$ such that $U=U(x, y)$; 
however, such a vector is not always unique. The claim follows from 
the fact that if
we denote the orbit through $(x, y)\in (\C^2\backslash 0)^k$ by $[x, y]$, then
the map $\Gamma_k\colon X_{\Sigma(k)}\rightarrow B(\C{[n-k+1, n+1]})$ defined by 
$\Gamma_k([x, y]):=U(x, y)$ is an equivariant diffeomorphism of complex 
manifolds for $1\leq k\leq n$ (see \cite{YC1}).     
\end{proof}
\begin{corollary}
The associated smooth fan $\Sigma(n)$ of the bounded flag manifold $B_n$ has
the generating set 
$G(\Sigma(n))=\{a_{1, 0}, a_{1, 1}, \ldots, a_{n, 0}, a_{n, 1}\}$, where
\begin{equation}
a_{j, \gamma}:=\begin{cases}
e_j& \text{if $\gamma=0$},\\
-e_1-\ldots -e_j& \text{if $\gamma=1$}
\end{cases}
\end{equation} 
for $1\leq j\leq n$.
\end{corollary}

%%%%%%%%%%%%%%%%%%%%%%%%%%%%%%%%%%%%%%%%%%%%%%%%%%%%%%%%%%%%%%%%%%%%%%%%%%%%
\subsection{Local Structures}\label{LS}
Even though we have displayed Bott towers as quotients, since we know
their associated smooth fans, we may also describe their local structures
in a toric setting, namely affine toric varieties forming each $X_{\Sigma(k)}$. 
In particular, we note that the collection of $k$-dimensional affine toric 
varieties provides an atlas for the manifold $X_{\Sigma(k)}$.
\begin{definition}
For a given $w=w_1\ldots w_k\in \SE^k$, we define
\begin{equation}
\mathbb{I}(j):=\{l\in [1, j)\colon w_l=1\}
\end{equation}   
for any $1\leq j\leq k$.
\end{definition}
\begin{definition}
Let a binary code $w=w_1\ldots w_k\in \SE^k$ be given. We define vectors in
$\R^k$ for each $1\leq j\leq k$ by
\begin{equation}
\upsilon_{j, w_j}:=\begin{cases}
\sum_{i\in \mathbb{I}(j)} (-1)^{w_j}\;b(\mathbb{O}^j_i)\;e_i 
+ (-1)^{w_j}\;e_j, & \text{if\;$\mathbb{I}(j)\neq \emptyset$},\\
(-1)^{w_j} e_j, & \text{otherwise}.
\end{cases}
\end{equation}
\end{definition}
\begin{proposition}
For any $w=w_1\ldots w_k\in \SE^k$, the dual cone $\Check{\sigma}_w$ of
$\sigma_w$ in $\Sigma(k)$ is given by
\begin{equation}
\Check{\sigma}_w= \RE \upsilon_{1, w_1} + \ldots + \RE \upsilon_{k, w_k}.
\end{equation}
\end{proposition}
\begin{proof}
We recall that $\Check{\sigma}_w=\{ x\in \R^k \colon \langle x, y \rangle \geq 0, 
\forall y\in \sigma_w\}$. If we set 
\begin{equation*}
d \sigma_w:=\RE \upsilon_{1, w_1} + \ldots + \RE \upsilon_{k, w_k},
\end{equation*}
then we need to show that $\Check{\sigma}_w= d \sigma_w$. 

So, let $x= x_1 \upsilon_{1, w_1} +\ldots + x_k \upsilon_{k, w_k}$ be any vector in
$d \sigma_w$. It follows that
\begin{equation*}
\langle x, y\rangle=
\sum_{i=1}^{k} \sum_{j=1}^{k} x_i y_j \langle \upsilon_{i, w_i}, a_{j, w_j}\rangle
\end{equation*}
for all $y=y_1 a_{1, i_1}+ \ldots + y_n a_{n, i_n}\in \sigma_w$.
By using \eqref{B2}, it can be verified that
\begin{displaymath}
\langle \upsilon_{i, w_i}, a_{j, w_j}\rangle=\begin{cases}
1,& \text{if $i=j$},\\
0,& \text{otherwise}
\end{cases}
\end{displaymath}
for any $1\leq i, j \leq k$. Therefore
$\langle x, y\rangle=\sum_{i=1}^{k} x_i y_i \geq 0$, providing
$x \in \Check{\sigma}_w$.

Conversely, let $x=(x_1, \ldots , x_k)$ be any vector in $\Check{\sigma}_w$ so that 
$\langle x, y\rangle\geq 0$ for all
$y\in \sigma_w$. If we define
\begin{displaymath}
u_j:=\begin{cases}
x_j& \text{if $w_j=0$},\\
b(j, j) x_j + \ldots + b(j, k) x_k& \text{if $w_j=1$},
\end{cases}
\end{displaymath}
then since $a_{j, w_j}\in \sigma_w$, it follows that
\begin{equation}
\langle x, a_{j, w_j} \rangle=u_j \geq 0
\end{equation}
for any $1\leq j\leq k$. We claim that 
$x= u_1 \upsilon_{1, w_1} +\ldots + u_k \upsilon_{k, w_k}$, providing 
$x \in d \sigma_w$. To see that we let
\begin{equation*}
v=( v_1, \ldots, v_k):= u_1 \upsilon_{1, w_1} +\ldots + u_k \upsilon_{k, w_k}.
\end{equation*}
If $w_j=0$ for some $1\leq j\leq k$, then $v_j= u_j= x_j$. On the other hand,
if $w_j=1$, then
\begin{align}\label{B4}
v_j=& -b(j, j)\;u_j + (-1)^{w_{j+1}} b(\mathbb{O}^{j+1}_j)\;u_{j+1} +
(-1)^{w_k} b(\mathbb{O}^k_j)\;u_k,\notag\\
=& x_j - \sum_{i=j+1}^k \Bigg(\sum_{l\in \mathbb{I}_j^i} 
b(\mathbb{O}^l_j)\;b(l, i) + b(j, i)\Bigg)\;x_i.
\end{align}
However, by \eqref{B3}, the coefficient of $x_i$ in \eqref{B4} for any
$j+1\leq i\leq k$ equals to zero; hence, $v_j=x_j$.
\end{proof}
\begin{corollary}\label{C3}
For any $w\in \SE^k$, the coordinate ring of the affine toric variety 
$X_{\Check{\sigma}_w}$ is given by $R_{\Check{\sigma}_w}=\C [\phi_1, \ldots,
\phi_k]\subset \C [z, z^{-1}]$, where $\phi_j:=z^{\upsilon_{j, w_j}}$ 
for each $1\leq j\leq k$.
\end{corollary}

%%%%%%%%%%%%%%%%%%%%%%%%%%%%%%%%%%%%%%%%%%%%%%%%%%%%%%%%%%%%%%%%%%%%%%%%%%
%%%%%%%%%%%%%%%%%%%%%%%%%%%%%%%%%%%%%%%%%%%%%%%%%%%%%%%%%%%%%%%%%%%%%%%%%%
\section{Classification Problem}\label{cla}

We devote this section to smooth toric varieties arising from crosspolytopes
and show that any such variety actually gives rise to a Bott tower. This may
be thought of as the generalization of Hirzebruch surfaces~\cite{TO}.
 
Throughout we will assume that $P^n$ is an $n$-crosspolytope in $\R^n$ 
with a set of vertices 
$V=\{x_1, \ldots, x_n, \widehat{x}_1, \ldots, \widehat{x}_n\}$, 
which spans a smooth projective fan $\Sigma^n=\Sigma(P^n)$. We therefore
identify the generating set of $\Sigma^n$ with the vertex set of $P^n$, i.e.,
$G(\Sigma^n)=V$.

We recall that any primitive collection of $\Sigma^n$ is of the form
$\PE=\{x_i, \widehat{x}_i\}$ for $i=1, \ldots, n$.
We then combine this fact with Proposition 3.2 of~\cite{Bat} to obtain the 
following result. 
\begin{corollary}\label{C2}
Let $\Sigma^n=\Sigma(P^n)$ be given as above. Then there exists a primitive
collection $\PE=\{x_i, \widehat{x}_i\}$ for some $1\leq i\leq n$ such that 
$x_i+\widehat{x}_i=0$.
\end{corollary}
Let us assume that $\PE^n=\{x_n, \widehat{x}_n\}$ is the primitive collection
of $\Sigma^n$ satisfying $x_n+\widehat{x}_n=0$, which exists by Corollary
\ref{C2}, and let $\sigma_n:=\RE x_n$ be the one-dimensional smooth cone of
 $\Sigma^n$ generated by $x_n$. We then consider the orthogonal projection
$p_n\colon \R^n\rightarrow \sigma^{\bot}_n=(\R x_n)^{\bot}$, and define
\begin{equation}
\Sigma^{n-1}:=\{p_n(\sigma)\colon \sigma\in
\Sigma^n\;\mathrm{and}\;\sigma_n\;\mathrm{is\;a\;face\;of}\;\sigma\}.
\end{equation}
In fact, $\Sigma^{n-1}$ is the quotient fan $\Sigma^n/ \sigma_n$
(see~\cite{GE}). Since $\Sigma^n$ is smooth, so is $\Sigma^{n-1}$. We may
proceed to define $\Sigma^k$ and $p_k$ for any $1\leq k <n-1$ in a similar
way in order to obtain a sequence of smooth fans:
\begin{equation}\label{E6}
\Sigma^n\xrightarrow{p_n} \Sigma^{n-1}\xrightarrow{p_{n-1}} \ldots 
\xrightarrow{p_3} \Sigma^2\xrightarrow{p_2} 
\Sigma^1 \xrightarrow{p_1} \{0\}.
\end{equation}
\begin{proposition}\label{P2}
For each $1\leq k\leq n-1$, the fan $\Sigma^k$ is a projective fan arising from
a crosspolytope.
\end{proposition}
\begin{proof}
We prove the claim when $k=n-1$ and a similar argument applies to the other
cases. Assume that $\PE^n=\{x_n, \widehat{x}_n \}$ is the primitive collection
for $\Sigma^n$ satisfying $x_n + \widehat{x}_n=0$. We observe that 
$\Sigma^{n-1}$ is a fan in $\sigma^{\bot}_n$ with the
generating set $G(\Sigma^{n-1})=\{y_1, \ldots, y_{n-1}, \widehat{y}_1, \ldots,
\widehat{y}_{n-1}\}$, where $y_i$ (respectively $\widehat{y}_i$) is the image
of $x_i$ (resp. $\widehat{x}_i$) under $p_n$ for each $1\leq i\leq
n-1$. By definition of a crosspolytope, each 
$P^i:=\mathrm{conv}\;\{x_1, \ldots, x_i, \widehat{x}_1, \ldots, 
\widehat{x}_i\}$ is an $i$-crosspolytope and satisfies
\begin{equation}
(\mathrm{relint}\;P^i)\cap (\mathrm{relint}\;[x_{i+1}, \widehat{x}_{i+1}])
=\{q_i\},
\end{equation} 
a single point, for $1\leq i\leq n-1$. By induction, we can obtain that
the set 

$R^i:=\mathrm{conv}\;\{y_1, \ldots, y_i, \widehat{y}_1, \ldots,
\widehat{y}_i\}$ is an $i$-crosspolytope in $\sigma^{\bot}_n$ and satisfies
\begin{equation}
(\mathrm{relint}\;R^i)\cap (\mathrm{relint}\;[y_{i+1}, \widehat{y}_{i+1}])
=\{p_n(q_i)\},
\end{equation} 
for each $1\leq i\leq n-2$; therefore, $\Sigma^{n-1}=\Sigma(R^{n-1})$. 
\end{proof}
Without any confusion we will assume that each fan $\Sigma^k$ is given by
$\Sigma^k=\Sigma(P^k)$, where $P^k$ is a $k$-crosspolytope with
$V (P^k)=\{x_1, \ldots , x_k, \widehat{x}_1, \ldots, \widehat{x}_k\}$, and
equivalently, $\PE^k=\{x_k, \widehat{x}_k \}$ is the
primitive collection of $\Sigma^k$ satisfying $x_k + \widehat{x}_k=0$
for each $1\leq k\leq n$.

We define $\Sigma^k_1$ and $\Sigma^k_2$ to 
be fans in $|\Sigma^k|=\bigcup_{\sigma \in \Sigma^k} \sigma\cong \R^k$ 
generated by the set of vectors $\{x_k, \widehat{x}_k \}$ and 
$\{x_1, \ldots , x_{k-1}, \widehat{x}_1, \ldots, \widehat{x}_{k-1}\}$ 
respectively. It can be easily seen that $\Sigma^k$ is the join of 
$\Sigma^k_1$ and $\Sigma^k_2$; that is, 
\begin{equation}\label{E7}
\Sigma^k=\Sigma^k_1 \cdot \Sigma^k_2.
\end{equation} 
Moreover, the map
$p_k \colon \Sigma^k \rightarrow \Sigma^{k-1}$ induces a bijection 
$p_k\arrowvert_{\Sigma^k_2} \colon \Sigma^k_2 \rightarrow \Sigma^{k-1}$;
that is, $\Sigma^{k-1}$ is the projective fan of $\Sigma^k$.    
Since the fan $\Sigma^k_1$ is unimodular equivalent to the fan of
$\C P^1$, we obtain the following. 
\begin{theorem}\label{T3}
Let $X^k$ denote the toric variety associated to the fan $\Sigma^k$ for each
$1\leq k\leq n$. Then, the induced toric morphism 
$\bar{p}_k \colon X^k \rightarrow X^{k-1}$ 
is a fiber bundle with fibers isomorphic to $\C P^1$.
\end{theorem}
Therefore, we can exhibit the associated smooth projective
toric variety $X^n$ as a sequence of smooth toric 
varieties of lower dimensions as a counterpart of \eqref{E6}:
\begin{equation}\label{E8}
X^n\xrightarrow{\bar{p}_n}
X^{n-1}\xrightarrow{\bar{p}_{n-1}}
\ldots \xrightarrow{\bar{p}_3}X^2\xrightarrow{\bar{p}_2}
X^1\xrightarrow{\bar{p}_1} *,
\end{equation}
where each map $\bar{p}_k$ is a fiber bundle with fibers isomorphic to 
$\C P^1$. However, this is very close to $\{X^k\colon k\leq n\}$ being a Bott
tower; in fact, we can say the following.
\begin{theorem}\label{T4}
If $\Sigma^n=\Sigma(P^n)$ is a smooth projective fan with $P^n$ being an
$n$-crosspolytope, then $\{X^k\colon k\leq n\}$ is a Bott tower.
\end{theorem}
\begin{proof}
We proceed by induction. If $k=1$, then $P^1$ is just a line segment so that
$X^1$ is isomorphic to $\C P^1$. So, assume for some
$1<k<n$, there exist line bundles $\xi_0, \ldots, \xi_{k-1}$ such that
$X^j= \C P(\C \oplus \xi_{j-1})$ for each $1\leq j\leq k$. 
Suppose $\Sigma^{k+1}$ is the associated fan to $X^{k+1}$ with the generating
set $G(\Sigma^{k+1})=\{x_1, \ldots, x_{k+1}, \widehat{x}_1, \ldots, 
\widehat{x}_{k+1}\}$ and the orthogonal projection 
$p_{k+1}\colon |\Sigma^{k+1}|\rightarrow (\R x_{k+1})^{\bot}$ such that
$\Sigma^k=\{p_{k+1}(\sigma)\colon \sigma\in
\Sigma^{k+1}\;\mathrm{and}\;\sigma_{k+1}\;\mathrm{is\;a\;face\;of}\;\sigma\}$
, where $\sigma_{k+1}=\RE x_{k+1}$. 
Let $\{u_1, \ldots, u_k\}$ be an orthonormal basis of 
$(\R x_{k+1})^{\bot}$. Then
the projection $p_{k+1}$ is given by 
$p_{k+1}(v)=\sum_{i=1}^k \langle v, u_i \rangle u_i$, and satisfies 
$v=p_{k+1}(v)+ (v- p_{k+1}(v))$ for each $v\in |\Sigma^{k+1}|$. We let  
$\{y_1, \ldots, y_k, \widehat{y}_1, \ldots, \widehat{y}_k\}$ be the generating
set of $\Sigma^k$ so that $y_i=p_{k+1}(x_i)$ and 
$\widehat{y}_i=p_{k+1}(\widehat{x}_i)$ for $1\leq i\leq k$. If we define a map
$h\colon |\Sigma^k|\rightarrow \R x_{k+1}\cong \R$ by 
\begin{equation}
h(y_i):=x_i- p_{k+1}(x_i)\;\;\mathrm{and}\;\;
h(\widehat{y}_i):=\widehat{x}_i- p_{k+1}(\widehat{x}_i),
\end{equation} 
the image of $\Sigma^k$ under the map $\Psi\colon |\Sigma^k|\rightarrow
|\Sigma^{k+1}|$ given by $y\mapsto (y, h(y))$ is exactly the fan obtained from
$\Sigma^{k+1}$ by removing its two one-dimensional cones $\RE x_{k+1}$ and
$\RE \widehat{x}_{k+1}$. On the other hand,
the map $h$ gives rise to an equivariant line bundle $L(h)$ 
over $X^k=\C P(\C\oplus \xi_{k-1})$ such that the space $\C P(\C \oplus L(h))$
is a toric variety with a fan being the join of 
$\Sigma_1:=\{\RE x_{k+1}, \RE \widehat{x}_{k+1}, \{0\}\}$ and
$\widetilde{\Sigma}:= \Psi(\Sigma^k)$. However, the fan $\Sigma_1\cdot
\widetilde{\Sigma}$ is in fact the fan $\Sigma^{k+1}$ associated to the toric
variety $X^{k+1}$. Therefore, we have $X^{k+1}=\C P(\C \oplus L(h))$, 
which completes the proof.  
\end{proof}
%%%%%%%%%%%%%%%%%%%%%%%%%%%%%%%%%%%%%%%%%%%%%%%%%%%%%%%%%%
%%%%%%%%%%%%%%%%%%%%%%%%%%%%%%%%%%%%%%%%%%%%%%%%%%%%%%%%%%%%%%
\section{Tangent Bundle of Bott Towers}\label{tb}
In this section, we investigate the tangent bundle of each $k$th stage $N_k$
of a Bott tower. Following~\cite{Ehl}, the tangent
bundle of any smooth toric variety has a natural decomposition into 
sum of line bundles, each of which arises from a one-dimensional 
cone of the associated fan. By Theorem \ref{T2}, we will not distinguish 
a Bott tower from its quotient description. One of the advantage of having 
such a description is that we may consider the action as an equivalence
relation, that is, any two vectors in $(\C^2\backslash 0)^k$ are equivalent
if they lie in the same orbit. Therefore, if $[x, y]$ denotes the equivalence
class of the vector $(x, y)\in (\C^2\backslash 0)^k$, then we may write
\begin{equation}\label{ER}
N_k=X_{\Sigma(k)}=\{[x, y]\colon (x, y)\in (\C^2\backslash 0)^k\}
\end{equation}
for each $k\leq n$.
\begin{definition}
Let a Bott tower $\{N_k\colon k\leq n\}$ of height $n$ be given and let
$\xi_0, \ldots, \xi_{n-1}$ be the associated line bundles. We then define 
$\lambda(k)$ to be the canonical line bundle over the projective space
$N_k=\C P(\C \oplus \xi_{k-1})$ for each $1\leq k\leq n$.   
\end{definition}
\begin{corollary}
For any $1\leq k\leq n$, there exists a line bundle $\lambda(k)^{\bot}$ over
$N_k$ satisfying
\begin{equation}\label{K8}
\lambda(k)\oplus \lambda(k)^{\bot}\cong \C \oplus \xi_{k-1}.
\end{equation}
\end{corollary}
\begin{proof}
Using the standard inner product in $\C \oplus \xi_{k-1}$, choose as a fiber
for $\lambda(k)^{\bot}$ over $L\in \C P(\C \oplus \xi_{k-1})$ the orthogonal
complement $L^{\bot}$ in $\C \oplus \xi_{k-1}$.
\end{proof}
Without any confusion we will continue to denote the line bundles over $N_k$
obtained by pulling back $\lambda_1$, $\ldots$, $\lambda_{k-1}$ along the
projections $p_j\colon N_j\rightarrow N_{j-1}$ for all $1\leq j\leq k$ and
$1\leq k\leq n$.
\begin{example}
Consider the Bott tower $\{B_k\colon k\leq n\}$ arising from bounded flag
manifolds and recall that $\eta_1$, $\ldots$, $\eta_{n-1}$ are the associated
line bundles. By the definition of the bundles $\eta_k$, we have 
$\lambda(k)=\eta_k$ and $\lambda(k)^{\bot}=\eta^{\bot}_k$
for each $1\leq k\leq n$.
\end{example}
Before proceeding further, we would like to give a detailed description of 
the bundles $\lambda(k)$ and $\lambda(k)^{\bot}$
over $N_k$. 
\begin{proposition}\label{K9}
For each $1\leq k\leq n$, the canonical bundle $\lambda(k)$ over $N_k$ is
given by
\begin{equation}
\lambda(k)\cong (\C^2\backslash 0)^k \times_{\C^k_{\times}} \C,
\end{equation} 
where the action of $\C^k_{\times}$ is defined by
\begin{equation*}
((x, y), v))\cdot t:=((x, y)\cdot t, t^{-1}_k v)
\end{equation*} 
for all $t=(t_1, \ldots, t_k)\in \C^k_{\times}$.
\end{proposition}
\begin{proof}
Since $N_k=\C P(\C \oplus\xi_{k-1})$, we may interpret an equivalence class
$[x, y]\in N_k$ as a line in $\C \oplus \xi_{k-1}$; hence, any vector 
$(x', y')$ in $[x, y]$ represents a point on this line. Therefore, the total
space of $\lambda(k)$ can be given as
\begin{equation*}
E(\lambda(k))=\{([x, y], (x', y'))\colon [x, y]\in N_k\;\;\mathrm{and}\;\;
(x', y')\in [x, y]\}.
\end{equation*} 
Now the map 
$g_k\colon (\C^2\backslash 0)^k \times_{\C^k_{\times}} \C\rightarrow
E(\lambda(k))$ defined by
\begin{equation*}
g_k([x, y; v]):=\begin{cases}
([x, y], (x_1, y_1; \ldots; x_{k-1}, y_{k-1}; x_k v^{-1}, y_k v^{-1})), 
& \text{if\;$v\neq 0$},\\
([x, y], (x_1, y_1; \ldots; x_k, y_k)),
& \text{if\;$v=0$}
\end{cases}
\end{equation*} 
provides the desired isomorphism.
\end{proof}
\begin{corollary}
Let a Bott tower $\{N_k\colon k\leq n\}$ be given with the integral sequence
$\mathbf{c}=(c_1, \ldots, c_n)$, and let $\xi_0, \ldots, \xi_{n-1}$ be the associated
line bundles. Then,
\begin{equation}\label{K11}
\xi_k\cong \lambda(1)^{-c(1, k+1)}\otimes \ldots \otimes
\lambda(k)^{-c(k, k+1)}
\end{equation} 
for all $1\leq k\leq n-1$.
\end{corollary} 
\begin{corollary}\label{C4}
For each $1\leq k \leq n$, the line bundle $\lambda(k)^{\bot}$ over $M_k$ can
be constructed as a quotient 
$(\C^2\backslash 0)^k\times_{\C^k_{\times}} \C$, where the action of
$\C^k_{\times}$ is given by 
\begin{equation*}
((x, y), v))\cdot t:=((x, y)\cdot t, t^{c(1, k)}_1\ldots 
t^{c(k-1, k)}_{k-1} t_k v)
\end{equation*}
\end{corollary}
\begin{proof}
When the isomorphism 
$\lambda(k)\oplus \lambda(k)^{\bot}\cong \C \oplus \xi_{k-1}$ is combined with
\eqref{K11}, we deduce that 
\begin{equation*}
\lambda(k)^{\bot}\cong \lambda(1)^{-c(1, k)}\otimes \ldots \otimes
\lambda(k-1)^{-c(k-1, k)}\otimes \lambda(k)^{-1},
\end{equation*}
from which the result follows.
\end{proof}
We are now ready to describe the tangent bundle of each $k$th stage $N_k$ of a
Bott tower. To do this, we recall the description of $N_k$ as a toric
variety $N_k\cong X_{\Sigma(k)}$, and the fact that the tangent bundle of
a toric variety has a decomposition into a sum of line bundles, 
each of which is constructed from a one-dimensional cone of the 
associated fan (see~\cite{Ehl}). Therefore, our next task is to 
identify these bundles for $X_{\Sigma(k)}$.

Let us denote the one-dimensional cones of $\Sigma(k)$ by
$\sigma^{\gamma}_j:=\RE a_{j, \gamma}$ for each $1\leq j\leq k$ 
and $\gamma=0, 1$, where
\begin{displaymath}
a_{j, \gamma}:=\begin{cases}
e_j& \text{if\;$\gamma=0$},\\
(0, 0, \ldots, 0, b(j, j), \ldots, b(j, k)) &\text{if\;$\gamma=1$}.
\end{cases}
\end{displaymath}
To simplify the notation, we abbreviate $\Sigma_{\sigma^{\gamma}_j}$ 
and $X_{\Sigma_{\sigma^{\gamma}_j}}$ to
$\Sigma^{\gamma}_j$ and $X^{\gamma}_j$ respectively, where
$\Sigma_{\sigma^{\gamma}_j}$ is the quotient fan of $\Sigma(k)$ by
$\sigma^{\gamma}_j$.
Then, it follows that the orthogonal complement 
of $\sigma^{\gamma}_j$ in $\R^k$ is $\R_{[k]\backslash \{j\}}$ 
if $\gamma=0$ and $V^k_j\times \R_{[j+1, k]}$ if $\gamma=1$, 
where $V^k_j$ is a $(j-1)$-dimensional subspace of $\R^k$. 
\begin{proposition}
Under the identification \eqref{ER}, the toric subvarieties 
$X^{\gamma}_j$ correspond to the codimension-one submanifolds 
$N^{\gamma}_j$ of $N_j$ defined by
\begin{itemize}
\item[(i)] 
$N^0_j:=\{[x, y]\in N_j\;|\; x_j=0\}$ if $\gamma=0$,
\item[(ii)]
$N^1_j:=\{[x, y]\in N_j\;|\; y_j=0\}$ if $\gamma=1$
\end{itemize}
for each $1\leq j\leq k$.
\end{proposition}
\begin{proof}
Let $\gamma=0$ and let $\sigma_w =\RE a_{1, w_1}+\ldots +\RE a_{k, w_k}$ be
any $k$-dimensional cone of $\Sigma(k)$ for some $w=w_1\ldots w_k\in \SE^k$
such that $\sigma^0_j$ is a face of it. Note that we then necessarily 
have $w_j=0$ in $w$. We see that any cone in $\Check{\Sigma}^0_j$ is of the form
\begin{displaymath}
\beta_w:=\Check{\sigma}_w \cap \R_{[k]\backslash \{j\}}=
\RE \upsilon_{1, w_1}+\ldots +\RE \upsilon_{j-1, w_{j-1}}+
\RE \upsilon_{j+1, w_{j+1}}+ \ldots +\RE \upsilon_{k, w_k}.
\end{displaymath}
Therefore, the associated affine toric subvariety $X_{\beta_w}$ of $X^0_j$ 
is embedded in $X_{\Check{\sigma}_w}=\{(\phi_1, \ldots, \phi_k)\}$
by assigning $\phi_j=0$.

On the other hand, $N^0_j$ being a codimension-one submanifold of $N_k$ 
has an open covering by the sets $X_{\Check{\sigma}_w}\cap N^0_j$ 
for any $w\in \SE^k$ for which $w_j=0$. Moreover, the intersection 
$X_{\Check{\sigma}_w}\cap N^0_j$ may be thought of the geometric 
realization of the affine variety $X_{\beta_w}$ providing the 
desired result. A similar argument applies to the case $\gamma=1$.
\end{proof}
\begin{example}
When we consider the Bott tower $\{B_k\colon k\leq n\}$ arising from bounded
flag manifolds, the submanifold $N^0_j$ of $B_n$ is
a copy of $B_{n-1}$ whose flags lie in $\C_{[n+1]\backslash \{j\}}$, and
similarly, $N^1_j$ is a copy of $B_{j-1}\times B_{n-j}$ whose flags lie in 
$\C_{[j]}\times \C_{[j+1, n+1]}$. 
\end{example}
\begin{proposition}
If we denote the bundle over $X_{\Sigma(k)}$ determined by the codimension-one 
subvarieties $X^{\gamma}_j$ by $\nu(\sigma^{\gamma}_j)$, there exist isomorphisms
of complex line bundles:
\begin{equation}
\nu(\sigma^0_j)\cong \bar{\lambda}(j)\quad\textrm{and}\quad  
\nu(\sigma^1_j)\cong \lambda(j)^{\bot},
\end{equation}
for each $1\leq j\leq k$, where $\bar{\lambda}(j)$ denotes the conjugate
of $\lambda(j)$.
\end{proposition}
\begin{proof}
This may be achieved directly by examining the corresponding transition
functions of these bundles (see \cite{Hir}). 
\end{proof}
\begin{theorem}
Let $\{N_k\colon k\leq n\}$ be a Bott tower of height $n$. Then
the tangent bundle $\tau(N_k)$ of each $k$th stage $N_k$ of the tower
satisfies 
\begin{equation}
\tau(N_k)\oplus \C^k \cong \bigoplus_{j=1}^k 
\bar{\lambda}(j)\oplus \lambda(j)^{\bot}.
\end{equation}
\end{theorem}
\begin{corollary}
When the bounded flag manifold $B_n$ is considered as a complex manifold,
there exists an isomorphism: 
\begin{equation}\label{EB}
\tau(B_n)\oplus \C^n \cong \bigoplus_{i=1}^n \bar{\eta}_j\oplus \eta_j^{\bot}.
\end{equation}
\end{corollary}
\begin{corollary}
If we denote the first Chern classes of line bundles $\lambda(j)$ and 
$\lambda(j)^{\bot}$ by $t_j$ and $z_j$ respectively, then the top Chern class
of $N_k$, the $k$th stage of the tower $\{N_k\colon k\leq n\}$ is given by
\begin{equation*}
c(N_k)=\prod_{j=1}^k (1-t_j+z_j).
\end{equation*}
\end{corollary}
We conclude our discussion with the cohomology ring of an arbitrary Bott 
tower by combining Theorem \ref{T2} and Danilov-Jurkiewicz theorem (see for
example~\cite{GE}). 
\begin{theorem}
Let $\{N_k\colon k\leq n\}$ be a Bott tower associated with the integral
sequence $\mathbf{c}=(c_1, \ldots, c_n)$. Then for each $k\leq n$, the cohomology ring 
$H^*(N_k; \Z)$ is isomorphic to the quotient $\Z[x_1, \ldots, x_k]/I_k$, where
$I_k$ denotes the ideal
\begin{equation*}
\bigg( x_j(x_j+ c(1, j)x_1+\ldots+c(j-1, j)x_{j-1})\colon 1\leq j\leq k \bigg).
\end{equation*}
\end{theorem}  
%%%%%%%%%%%%%%%%%%%%%%%%%%%%%%%%%%%%%%%%%%%%%%%%%%%%%%%%%%%%%%%%%%%%%%
%%%%%%%%%%%%%%%%%%%%%%%%%%%%%%%%%%%%%%%%%%%%%%%%%%%%%%%%%%%%%%%%%%%%%%%

\end{document}